\input amstex
\documentstyle {amsppt}
\def\picture #1 by #2 (#3 scaled #4){\vskip 6pt\centerline{
 \vbox to #2{
    \hrule width #1 height 0pt depth 0pt

    \vfill
    \special{picture #3 scaled #4} 
}}
\vskip 6pt
}

\advance\voffset  by -1.0cm
\NoBlackBoxes
\magnification=\magstep1
\hsize=17truecm
\vsize=24.2truecm
\voffset=0.5truecm
\document

\define \Ga {\Gamma}
\define \la {\lambda}
\define \pa {\partial}
\define \bR {\Bbb R}
\define \bP {\Bbb P}
\define \bC {\Bbb C}
\define \al {\alpha}
\define \ga {\gamma}
\def \Ga {\Gamma}
\define \CC {\frak {Con}}
\define \ND {\frak {ND}}
\define \GEN {\frak {GEN}}
\define \stp {S^1\to \bP^n}
\topmatter
\title
Discriminants of convex curves are homeomorphic
\endtitle

\author B.~Shapiro
\endauthor
\affil
 Department of Mathematics, University of Stockholm, S-10691, Sweden,
{\tt shapiro\@matematik.su.se}\\
\endaffil
\rightheadtext {Discriminants of convex curves}

\abstract For a given real generic curve $\ga: S^1\to \bP^n$ let $D_\ga$
denote the ruled hypersurface in $\bP^n$ consisting of all osculating
subspaces to $\ga$ of codimension 2.
In this short note we show that for any two convex real projective curves
$\ga_1:S^1\to\bP^n$ and $\ga_2:S^1\to\bP^n$  the pairs $(\bP^n,D_{\ga_1})$ and
$(\bP^n,D_{\ga_2})$ are homeomorphic. \endabstract

\subjclass {Primary 14H50}\endsubjclass
\keywords{convex curves, discriminants}
\endkeywords

\endtopmatter

\heading  \S 0. Preliminaries and results\endheading

{\smc Definition.} A smooth curve $\ga:S^1\to \bP^n$ is called {\it
nondegenerate
or locally convex}  if the local multiplicity
of its intersection with any hyperplane does not exceed $n$, i.e. in local
terms
$\ga^\prime(t),...,\ga^{(n)}(t)$ are linearly independent at every $t$ or its
osculating complete flag is well-defined at every point. A curve
$\ga:S^1\to \bP^n$
is called {\it convex}  if the total multiplicity
of its intersection with any hyperplane does not exceed $n$.

The set $\CC_n$ of all convex curves in $\bP^n$ forms 1 connected component of
the space $\ND_n$ of
all nondegenerate curves if $n$ is even and 2 connected components (since
the osculating
frame orients $\bP^{2k+1}$) if $n=2k+1$, see  \cite {MSh}. Different
results about
convex curves show that they have the most
simple  properties among all curves. In this paper we prove one more result
of the same
nature.

{\smc Definition.\,} A curve $\ga:S^1\to \bP^n$ is called {\it generic}  if
at every point
 $\ga(t),\,t\in S^1$ one has a well-defined osculating subspace of
codimension 2, i.e.
in local terms $\ga^\prime(t),...,\ga^{(n-1)}(t)$ are linearly independent
at every $t$.

Note that any smooth curve $\ga:\stp$  can be made
generic by a small smooth deformation of the map. The space  $\ND_n$
of all nondegenerate curves is the union of some connected components in
the space
$\GEN_n$ of all generic curves. (The number of
connected components in $\ND_n$
 equals 10 for odd $n\geqslant 3$ and equals 3 for even $n\geqslant 2$, see
\cite {MSh}.)

{\smc Definition.} Given a generic $\ga:S^1\to \bP^n$ we call by its {\it
standard discriminant
$D_\ga\subset \bP^n$} the hypersurface consisting of all codimension 2
osculating
subspaces to $\ga$.

In many cases (algebraic, analytic etc) the assumption of genericity in the
definition of discriminant can be omitted.

The following proposition answers the question posed by V.Arnold in \cite
{Ar2}, p.37.

\proclaim {Main proposition} a) For any 2 convex curves
$\ga_1: S^1\to \bP^n$ and $\ga_2: S^1\to \bP^n$  the pairs
$(\bP^n,D_{\ga_1})$ and $(\bP^n,D_{\ga_2})$ are homeomorphic.

b) For any
convex curve $\ga$ the complement $\bP^n\setminus D_\ga$ consists of
$[\frac n2] +1$ components. All components are contractible to $S^1$ for $n$
even and all but one are contractible to $S^1$ for $n$ odd. The
remaining component is a cell.
\endproclaim

Now we want to place this result into a more general context of associated
discriminants
in the spaces of (in)complete flags.

{\smc Notation.} Let $F_{n+1}$ denote the space of all complete flags in
$\bP^n$ (or,
equivalently, in $\bR^{n+1}$). Given a nondegenerate curve
$\ga:S^1\to\bP^n$ one can consider
its {\it associated curve} $\tilde\ga: S^1\to F_{n+1}$ where $\tilde
\ga(t)$ is the complete
osculating flag to $\ga\,$ at $\ga(t)$. Note that any associated curve
$\tilde\ga: S^1\to F_{n+1}$
is tangent to the special distribution of $n$-dimensional cones in
$F_{n+1}$ and any integral
curve of this distribution is the associated curve of some nondegenerate
projective curve,
see e.g. \cite {Sh1}.

Given a complete flag $f\in F_{n+1}$ and some space
$\frak G=SL_{n+1}/P\,$ of (in)complete flags  where $P$ is some parabolic
subgroup one gets the
Schubert cell decomposition $\frak {Sch}_f$ of $\frak G$ as follows. Each
cell of $\frak {Sch}_f$
consists of all flags in $\frak G$ subspaces of which have a given set of
dimensions of
intersections with the subspaces of $f$. Let $\frak O_f$ denote the union
of all cells in
$\frak {Sch}_f$ which have codimension at least 2. (Obviously, $\text
{codim  } \frak O_f=2$.)

{\smc Examples.} 1) If $\frak G$ equals $\bP^n$ then $\frak O_f$ is the
subspace of $f$ of
codimension 2.

2) If $\frak G=F_3$ then $\frak O_f$ consists of 2 copies of $\bP^1$
intersecting
at $f$. The first $\bP^1$ is the set of flags on $\bP^2$ with the same
point as that of $f$
and the second $\bP^1$ consists of all flags with the same line as that of $f$.

{\smc Definition.} For a given curve $c:S^1\to F_{n+1}$ and a space $\frak
G=SL_{n+1}/P\,$
of (in)complete flags we denote by its {\it $\frak
G$-discriminant\,} $\frak{GD}_c$ the union
$\bigcup_{t\in S^1} \frak O_{c(t)}\subset \frak G$. (If $c$ is not a
constant map then $\frak {GD}_c$ is a
hypersurface in $\frak G$.)

Note that the standard discriminant $D_\ga$ of a nondegenerate curve
$\ga:S^1\to\bP^n$
can be considered as the $\frak G$-discriminant for $\frak G=\bP^n$.

{\smc Definition.} Two nondegenerate curves $\ga_1: S^1\to \bP^n$ and
$\ga_2: S^1\to \bP^n$
are called {\it $\frak G$-equivalent\,} if the pairs $(\frak G,
\frak{GD}_{\tilde \ga_1})$ and
$(\frak G, \frak{GD}_{\tilde \ga_2})$ are homeomorphic. (Recall that
$\tilde \ga$ denotes the
associated curve of $\ga$.)

{\smc Remark.} The notion of $\frak G$-equivalence of nondegenerate curves
is intrinsically related with the
qualitative theory of linear ODE since each nondegenerate curve in $\bP^n$
can be represented as the
projectivization of the fundamental solution of some linear ODE of order
$n+1$, see \cite {Sh2}.
The problem of enumeration of $\frak G$-equivalent generic curves
 is apparently a very interesting and difficult question even for
$n=2$.

The following conjecture is formulated in \cite {Sh2}.

\proclaim {Conjecture 1} Any 2 convex curves
$\ga_1: S^1\to \bP^n$ and $\ga_2: S^1\to \bP^n$ are $\frak G$-equivalent
for any $\frak G=SL_{n+1}/P\,$.
\endproclaim

Note that it suffices to prove this conjecture in the case of the space
$F_{n+1}$ of complete flags, i.e. to the case $P=B$ where $B$ is the Borel
subgroup of uppertriangular matrices.

\proclaim {Conjecture 2} For any pair of convex curves
$\ga_1: S^1\to \bP^n$ and $\ga_2: S^1\to \bP^n$  the pairs
$(\bP^n,D_{\ga_1})$ and $(\bP^n,D_{\ga_2})$ are diffeomorphic.
\endproclaim

The main motivation of this paper was an attempt to formalize the  idea
that any 2 convex curves
as qualitatively equivalent in any natural sense.
It is difficult to overestimate the role of my visit to the Max-Planck
Institute during the
summer 1996 where the main bulk of this project was carried out.
Stimulating discussions with S.~Tabachnikov and Vl.~Lin are highly
appreciated.

\heading  \S 1. Proofs. \endheading

\subheading {Some generalities on convex curves}

{\smc Definition.}  For any $t\in S^1$ and $1\le k\le n-1$ let $L^k_t$
denote the osculating subspace to $\ga$ at $\ga(t)$ of dimension $k$.

\proclaim {1.1. Theorem (criterion of convexity)} A curve $\ga:\stp$ is convex
if and only if
for any $r$-tuple of positive integers $k_1,..., k_r$ such that $\sum k_i=n$ and
any $r$-tuple of pairwise different moments $t_1,...,t_r$ the intersection
$L_{t_1}^{n-k_1}\cap ... \cap L_{t_r}^{n-k_r}$ is a point.
\endproclaim

\demo {Proof} In order to save the space we refer the interested reader to
 \cite {Co},\cite {Ma} and further references.

\qed
\enddemo

{\smc Definition.} Given a nondegenerate curve $\ga:\stp$ we call by its dual
 $\ga:S^1\to (\bP^n)^*$ the curve consisting of all osculating hyperplanes
to $\ga$.

{\smc Remark.} If $\ga$ is convex then $\ga^*$ is also convex, see \cite {Ar1},
\cite {Ar2}.

{\smc Notation.} If a point $p$ lies on the osculating hyperplane $H_\tau$
to $\ga$
we say that {\it the order of tangency
$\sharp_p(\ga(\tau))$} of $p$ at $\ga(\tau)$  equals to $i$
if $p$ belongs to the osculating subspaces at $\ga(\tau)$ of
the codimension at most $i$.  ( For example, for every point $p$ on a line $l$
tangent to a circle $c$ at $c(1)$ on $\bP^2$ except for the tangency point
$c(1)$ 
one has $\sharp_p(c(1))=1$. On the other side, $\sharp_{c(1)}(c(1))=2$.)

 Given a nondegenerate $\ga:\stp$ and a point $p\in\bP^n$
we call
by {\it the number of roots $\sharp_p(\ga)$} of $p$ the sum of the
orders of tangency $\sharp_p(\ga(t_i))$ taken over all osculating
hyperplanes $H_{t_i}$ through $p$.

(The term 'number of roots' comes from the example when $\ga$ is the rational
normal curve in $\bP^n$. In this case all points in $\bP^n\,$ can be
interpreted
as homogeneous polynomials in 2 variables of degree $n$ with real coefficients
(considered up to a constant factor)
and $\ga$ is the family of polynomials of the form
$(ax_1+bx_2)^n,\,a^2+b^2\neq 0$.
In this  situation
$\sharp_p(\ga)$ coincides with the total number of real roots of such a
 polynomial on $\bP^1$ counted with multiplicies.)

\proclaim {Observation} The number of roots $\sharp_p(\ga)$ coincides with the
total multiplicity (i.e. sum of local multiplicities) of the intersection of
$\tilde H_p\,$ with $\ga^*$. Here $\tilde H_p$ denotes the hyperplane in
$(\bP^n)^*$ corresponding to the point $p\in \bP^n$. \endproclaim

\proclaim {1.2. Corollary} A curve $\ga$ is convex if and only if for any
$p\in \bP^n$
one has $\sharp_p(\ga)\le n$. \endproclaim

\subheading {Projection}

Given a convex curve $\ga:\stp$ and its osculating hyperplane $H_\tau$ at
the point
$\ga(\tau)$ let us denote by $\ga^\tau:S^1\to H_\tau$ the curve   obtained by
projection of $\ga\,$ onto $H_\tau$ along the pencil of tangent lines to $\ga$,
i.e. for any $t\in S^1$ one has $\ga^\tau(t)=H_\tau\cap l_t$ where $l_t$ is
the tangent line to $\ga$
at $\ga(t)$.

\proclaim {1.3. Lemma} For any $\tau\in S^1$ the curve $\ga^\tau$ is a convex
curve in $H_\tau$. Osculating hyperplanes to $\ga^\tau$ and its discriminant
$D_{\ga^\tau}$ are obtained by intersection of the osculating  hyperplanes and
$D_\ga\,$ with $H_\tau$.
\endproclaim

\demo {Proof} The argument splits into 2 principal parts. At first we show that
$\ga^\tau$ is nondegenerate, i.e.
$(\ga^\tau)^\prime(t),...,(\ga^\tau)^{(n-1)}(t)$ are
linearly independent at any $t\in S^1$. Then we prove that $\ga^\tau$ is
convex, i.e. its total
multiplicity of intersection with any hyperplane in $H_\tau$ does not exceed
$n-1$. Observe that $\ga$ has the only intersection point with $H_\tau$, namely
$\ga(\tau)$. Assume first that $t\neq \tau$. In this case one can choose
a system of affine coordinates $x_1,...,x_n$ in $\bP^n$ such that  $H_\tau$
coincides with the hyperplane $\{x_n=0\}$; $\ga(t)$ is the point with
coordinates
$(0,...,0,1)$ and the tangent line $l_t$ to $\ga\,$ at $\ga(t)$ is the
$x_n$-axis.
In these coordinates
the curve $\ga^\tau$ has the form
$$\ga^\tau(t)=\ga(t)-\frac{\ga_n(t)}{\ga_n^\prime(t)}\ga^\prime(t)$$
where $\ga_n$ is the last coordinate of $\ga$. (Under our assumptions
$\ga_n^\prime(t)\neq 0$.)
Therefore
 $$(\ga^\tau)^{(i)}(t)=(-1)^i\frac {\ga_n(t)}{\ga_n^\prime(t)}\ga^{(i)}(t)+...$$
where $...$ denotes the terms containing derivatives of $\ga$ of order
lower than $i$.
By the above assumptions $\frac {\ga_n(t)}{\ga_n^\prime(t)}\neq 0$ and since
$\ga^\prime(t),...,\ga^{(n)}(t)$ are linearly independent one gets that the
derivatives
$(\ga^\tau)^{(i)}(t),\;i=1,...,n-1$ are linearly independent as well.

The alternative geometric argument is as follows. Since the point $\ga(t)$
does not lie on $H^\tau$ one has that the osculating complete flag $f(t)$ to
$\ga$ at $\ga(t)$ is transversal to $H_\tau$ and the same holds for all
$t^\prime$
close to $t$. Therefore the complete flags obtained by intersection of
$f(t^\prime)$
with $H_\tau$ are well defined. But in their turn these flags coincide with
the osculating
flags to $\ga^\tau$ which are therefore well-defined in some neighborhood
of $t$.

It is left to show that $\ga^\tau$ is nondegenerate at $t=\tau$. This follows
from the local calculation given below. In this case we can choose a system of
coordinates such that in the neighborhood of $\tau$ (we assume $\tau=0$) the
curve $\ga$ has the form
$$x_1=t+...,\,x_2=t^2+...,\,...,x_n=t^n+...$$
The osculating hyperplane $H^0$ at $\tau=0$ is given by $\{x_n=0\}$.
The projected curve $\ga^0(t)$ is given by
$$\align \ga^0(t) &=\ga(t)-\frac {\ga_n(t)}{\ga_n^\prime(t)}\ga^\prime =
(t+...,t^2+...,...,t^n+...)-\frac {t^n+...}{nt^{n-1}+...}(1+...,2t+...,...,
nt^{(n-1)}+...)\\
&=\frac 1n((n-1)t+...,(n-2)t^2+...,...,t^{(n-1)}+...,0)\endalign$$
which shows that $\ga^0$ is nondegenerate at $t=0$.

Now we show that $\ga^\tau$ is convex. By Corollary 1.2. one has to prove that
for any $p\in H_\tau$ the number of roots $\sharp_p(\ga^\tau)$ is less or
equal $n-1$.
This follows from the equality
$$\sharp_p(\ga^\tau)+1=\sharp_p(\ga)$$
which together with convexity of $\ga$ gives the required result.
Indeed, assume that some $p\in H_\tau$ lies in the intersection
$L^{k_1}_{t_1}(\ga)\cap
 L^{k_2}_{t_2}(\ga)\cap ...\cap L^{k_r}_{t_r}(\ga)$ where $t_1=\tau$. Since
each subspace
$L^{k_i}_{t_i}$ for $i\neq 1$ is transversal to $H_\tau$ (see criterion of
convexity)
one has that $p\,$ lies in the intersection $L^{k_1}_{t_1}(\ga^\tau)\cap
 L^{k_2-1}_{t_2}(\ga^\tau)\cap ...\cap L^{k_r-1}_{t_r}(\ga^\tau)$.
Therefore, by definition
of the number of roots, one  gets the above equality.

\qed
\enddemo

 For any $k$-tuple of moments $(t_1,...,t_k),\,t_i\in S^1$ let $H_{t_1}\cap
... \cap H_{t_k}$
denote the intersection of the osculating hyperplanes
$H_{t_i},\,i=1,...,k$. In what follows
we use the following convention.
 If some of the
 moments $t_{j_1},t_{j_2},...,t_{j_r}$ coincide we define the intersection
$H_{t_{j_1}}\cap H_{t_{j_2}}\cap ...\cap H_{t_{j_r}}$ as the osculating
subspace to $\ga$
at $\ga(t_{j_1})$ of codimension $r$. Under this convention one has that
$H_{t_1}\cap ... \cap H_{t_k}$
always has codimension $k$, see 1.1.

\proclaim {1.4. Corollary} The projection  $\ga^{t_1,...,t_k}\,$ of $\ga$
onto any
intersection of osculating hyperplanes $H_{t_1}\cap ... \cap H_{t_k}$  by a
pencil of
$k$-dimensional osculating subspaces to
$\ga$  is a convex curve. For any point $p\in H_{t_1}\cap ...
\cap H_{t_k}$ one has
$$\sharp_p(\ga^{t_1,...,t_k})+k=\sharp_p(\ga).$$
\endproclaim

\demo {Proof} Apply the above lemma several times.

\qed
\enddemo

\subheading {Elliptic hull of $\ga$ and root filtration of $\bP^n$}

{\smc Definition.} For a convex $\ga:\stp$ we denote by {\it its elliptic hull
$Ell_\ga$}  the set of all $p\in \bP^n$ with
$$\cases  \sharp_p(\ga)=0,\text { if $n$ is even}\\
   \sharp_p(\ga)=1, \text { if $n$ is odd.}  \endcases $$

\proclaim {1.5. Lemma} a) If $n\,$ is even then $Ell_\ga\,$ is a nonempty
convex set
in some affine chart of $\bP^n$, comp. \cite {ShS};

b) If $n$ is odd then $Ell_\ga\,$ is a disjoint union of $\bigcup_{\tau\in S^1}
Ell_{\ga^\tau}$ and, therefore, is fibered over $\ga$ with a contractible
fiber.
\endproclaim

\demo {Proof} a) Note that if $\ga:S^1\to \bP^{2k}$ is convex then $\ga$ lies in
some affine chart in $\bP^{2k}$. Indeed, take some osculating hyperplane
$H_\tau$.
The curve $\ga$ is tangent to $H_\tau$ only at $\ga(\tau)$ with the multiplicity
$2k$. Locally $\ga$ lies on one side w.r.t. $H_\tau$. Therefore, one can
make a small
shift of $H_\tau$
in order to get rid of the intersection points with $\ga$ near $\ga(\tau)$.
But no new
intersection can appear for a sufficiently small shift since the only
intersection point of
$\ga$ and $H_\tau$ is $\ga(\tau)$.

Assume now that $\bR^{2k}\subset \bP^{2k}$ is
the affine
chart containing $\ga$. We claim that $Ell_\ga$ coincides with the intersection
$\bigcap_{\tau\in S^1} Half_\tau$. Here $Half_\tau$ is the open halfspace in
$\bR^{2k}$ containing $\ga$ and bounded by the osculating hyperplane $H_\tau$.
First of all, $\bigcap_{\tau\in S^1} Half_\tau$ is nonempty since it is an open
 convex set containing the interior of the convex hull of $\ga$ in $\bR^{2k}$.
Then $\bigcap_{\tau\in S^1} Half_\tau$ is
contained in $Ell_\ga$. Indeed, every hyperplane through a point $p\in
\bigcap_{\tau\in S^1} Half_\tau$ is transversal to any osculating
hyperplane $H_\tau$ since $p\notin H_\tau$. On the other side,
$Ell_\ga\subseteq \bigcap_{\tau\in S^1} Half_\tau$. Indeed, for every
$p\notin Half_\tau ,\,$
  $\tau\in S^1$ there exists a hyperplane $L_p$ through $p$ not intersecting
$\ga$ at all. Take the affine chart $\bP^n\setminus L_p$ containing $\ga$ 
and some pencil $\Cal L_p$ of 'parallel' hyperplanes through $p$. Since
$\ga$ is a closed
curve in $\bP^n\setminus L_p$ one gets that some hyperplane in $\Cal L_p$
does not intersect $\ga$. Therefore, there exists a hyperplane in $\Cal L_p$
tangent to $\ga$ at some $\ga(t_p)$. But this exactly means that the osculating
hyperplane $H_{t_p}$ contains $p$.

b) Take a 1-parameter family of osculating hyperplanes.
According to the proof of lemma 1.3. for any $\tau\in S^1$
the curve $\ga^\tau$ is convex in $H_\tau$ and one has
$\sharp_p(\ga^\tau)+1=\sharp_p(\ga)$. Therefore
the elliptic domain $Ell_{\ga^\tau}$ of every  curve $\ga^\tau$ belongs to
$Ell_\ga$, i.e. $\bigcup_{\tau\in S^1}Ell_{\ga^\tau}\subset Ell_\ga$. (Note
that
the union $\bigcup_{\tau\in S^1}Ell_{\ga^\tau}$ is disjoint.) Conversely,
 by definition, for odd $n$  every point $p$ in $Ell_\ga$  has
exactly one tangent hyperplane to $\ga$ and thus $p$ belongs  exactly to
one osculating
$H_\tau$. By the equality $\sharp_p(\ga^\tau)+1=\sharp_p(\ga)$ the point
$p$ lies in the
elliptic hull of $\ga^\tau$. Moreover, by the first part of this proof,
 $Ell_{\ga^\tau}$ is a convex domain in $H_\tau$ and, therefore, is
contractible which
gives the necessary result.

\qed
\enddemo

{\smc Definition.} By {\it the root filtration} of $\bP^n$ w.r.t. a convex
curve $\ga:\stp$
$$\bP_0(\ga)\subset ... \subset \bP_{[\frac n 2 ]}(\ga)=\bP^n$$
we denote  the filtration where each $\bP_i(\ga)$ consists of all $p\in \bP^n$
for which the number of roots $\sharp_p(\ga)\,$ is greater or equal $n-2i$.

Let $\Cal T^j=(S^1)^j$ denote the $j$-dimensional torus and let
$\Cal T^j/\frak S_j$ be its quotient mod the natural action of the
symmetric group
$\frak S_j$ by permutation of copies of $S^1$.

\proclaim {1.6. Lemma} a) For any $n$ and  $0\le i\le [\frac n 2]$ the set
$\bP_i(\ga)\setminus \bP_{i-1}(\ga)$ is  naturally fibered over
$\Cal T^{n-2i}/\frak S_{n-2i}$
with a contractible fiber. (For $n=2k$ and $i=k$ the set
$\bP_{k}(\ga)\setminus \bP_{k-1}(\ga)$ is contractible, see 1.5.a);

b) This fibration is trivial.

\endproclaim

\demo {Proof} a) Every point $p\in \bP_i(\ga)\setminus \bP_{i-1}(\ga)$ can be
described as follows.
There exists and unique $(n-2i)$-tuple of osculating hyperplanes $H_{t_1},...,
H_{t_{n-2i}}$ to $\ga$ (with probably coinciding moments $t_1,...,t_{n-2i}$
in which case
we use the same convention as above) such that $p$ belongs to the intersection
$H_{t_1}\cap H_{t_2}\cap ...\cap H_{t_{n-2i}}$ and, moreover, lies in the
elliptic hull of the curve $\ga^{t_1,...,t_{n-2i}}$. (Here
$\ga^{t_1,...,t_{n-2i}}$
is the projection of $\ga$ onto $H_{t_1}\cap H_{t_2}\cap ...\cap
H_{t_{n-2i}}$
by the pencil of osculating subspaces of dimension $n-2i$.) Indeed, we have
that $\sharp_p
(\ga^{t_1,...,t_{n-2i}})+2i=\sharp_p(\ga)$, see 1.2. Therefore $p$ must lie
in the
elliptic hull
of $\ga^{t_1,...,t_{n-2i}}$. On the other side, any intersection
$H_{t_1}\cap H_{t_2}\cap ...\cap H_{t_{n-2i}}$ has
codimension $n-2i$,
see 1.1. and any curve $\ga^{t_1,...,t_{n-2i}}$ is convex.
Therefore applying 1.5. we get the necessary result.

b) The fibration of elliptic components $Ell_{\ga^{t_1,...,t_{n-2i}}}$ of
the curves $\ga^{t_1,...,t_{n-2i}}$ over the set of moments
$(t_1,...,t_{n-2i})\in
\Cal T^{n-2i}/\frak S_{n-2i}$ depends continuously on $\ga\in \frak
{Con}_n$. Since
$\frak {Con}_n$ consists of 1 connected component (up to orientation for
odd $n$)
it suffices to show that the fibration sending $Ell_{\ga^{t_1,...,t_{n-2i}}}$
to $(t_1,...,t_{n-2i})$
is trivial for some $\ga\in \frak {Con}_n$.

The simplest example showing triviality is the case when $\ga$ is the
rational normal curve.
Indeed, in this case the space under consideration is the fibration of the
space
$\Pi_n(i)$ all homogeneous
forms of degree $n$ in 2 variables (up to a scalar multiple) which have
exactly $n-2i$
real zeros counted with multiplicities over the space $\Cal T^{n-2i}/\frak
S_{n-2i}$
of their real zeros. But $\Pi_n(i)$  has the obvious structure of the
product of the
space of  degree $n-2i$ polynomials
with all real zeros  (considered up to a scalar multiple) and the space of
degree $2i$
 polynomials
 with no real zeros (considered up to a scalar multiple). This shows that
the fibration
$\Pi_n(i)\to \Cal T^{n-2i}/\frak S_{n-2i}$ is trivial.

\qed
\enddemo

\subheading {Proof of the main proposition}

a) We will construct the homeomorphism of pairs $(\bP^n,D_{\ga_1})$ and
$(\bP^n,D_{\ga_2})$
in $[\frac n 2] +1$ steps. On the $i$th step, $i=0,[\frac n 2]$ we obtain
the partial
homeomorphism $h_i$ of the terms  $\bP_i(\ga_1)$ and $\bP_i(\ga_2)$ of the
above filtration.

The initial step. We construct the homeomorphism $h_0:\bP_0(\ga_1)\to
\bP_0(\ga_2)$.
Indeed, each of $\bP_0(\ga_1)$ and  $\bP_0(\ga_2)$ is homeomorphic to $\Cal
T^n/\frak S_n$
as follows. Every element in $\Cal T^n/\frak S_n$  is a pair
$(t_1,...,t_r)\in (\Cal T^r\setminus Diag)/\frak S_r,r\le n$ and
$(k_1,...,k_r),\,\sum k_i=n$. We map such a pair
$(t_1,...,t_r),(k_1,...,k_r)$ onto the intersection point
$L^{n-k_1}_{t_1}(\ga_j)\cap ...
\cap L_{t_r}^{n-k_r}(\ga_j), j=1,2$. This identification provides the
homeomorphism
$h_0:\bP_0(\ga_1)\to \bP_0(\ga_2)$ by 1.1.

The typical step. Each point in $\bP_i(\ga_j)\setminus \bP_{i-1}(\ga_j)$
lies in
the eliptic hull of the unique curve $\ga^{t_1,...,t_{n-2i}}\subset
H_{t_1}\cap ...
\cap H_{t_{n-2i}}$, i.e. the set of (nonnecessarily pairwise different) moments
$(t_1,...,t_{n-2i})\in \Cal T^{n-2i}/\frak S_{n-2i}$ is uniquely defined.
For each individual intersection
$ H_{t_1}\cap ... \cap H_{t_{n-2i}}$  the homeomorphism $h_{i-1}$ is already
defined on the complement to the elliptic hulls of the curves
$\ga_1^{t_1,...,t_{n-2i}}$ and $\ga_2^{t_1,...,t_{n-2i}}$. Since the
elliptic hulls are
convex domains and the fibrations of the elliptic hulls over
$\Cal{T}^{n-2i}/\frak S_{n-2i}$ 
are trivial we can extend  $h_{i-1}$ fiberwise to $h_i$ by identifying the
points  of the
elliptic hull of $\ga_1^{t_1,...,t_{n-2i}}$ with points of the elliptic hull of
$\ga_2^{t_1,...,t_{n-2i}}$  for all tuples $(t_1,...,t_{n-2i})\in
\Cal{T}^{n-2i}/\frak S_{n-2i}$

b) The corresponding component $Comp_i$ of $\bP^n\setminus D_\ga$ contained in
$\bP_i\setminus \bP_{i-1}$ is fibered over  $(\Cal T^{n-2i}\setminus
Diag)/\frak S_{n-2i}$
with the contractible fiber.  Since $(\Cal T^{n-2i}\setminus Diag)/\frak
S_{n-2i}$ is
contractible to $S^1$ for any $n-2i>0$ one gets that  $Comp_i$ is
contractible to $S^1$
for all $n$ and $i\le [\frac n 2]$ except for $Ell_\ga$ for even $n$
which is contractible to a point.

\qed

\enddocument